\pdfoutput=1 
\documentclass[12pt]{article}
\usepackage{amsmath,amsfonts,amssymb,url,color,epsfig,graphics}
\usepackage[latin1]{inputenc}

\newcommand{\R}{{\mathbb{R}}}

\newcommand{\C}{{\mathbb{C}}}
\newcommand{\Z}{{\mathbb{Z}}}
\newcommand{\N}{{\mathbb{N}}}

\def\ha{\frac{1}{2}}

\def\pa{\partial}
\def\ra{\rightarrow}

\def\nghbd{neighbourhood~}

\newtheorem{defi}{Definition}[section]

\newtheorem{lemm}{Lemma}[section]
\newtheorem{prop}{Proposition}[section]
\newtheorem{rem}{Remark}[section]

\newtheorem{theo}{Theorem}[section]

\newtheorem{conj}{Conjecture}[section]

\begin{document}

\title{On essential-selfadjointness of differential operators on
closed manifolds}

\author{Yves  Colin de Verdi\`ere\footnote{Institut Fourier, Universit\'e Grenoble-Alpes,
 Unit{\'e} mixte
 de recherche CNRS-UGA 5582,
 BP 74, 38402-Saint Martin d'H\`eres Cedex (France);
{\color{blue} {\tt yves.colin-de-verdiere@univ-grenoble-alpes.fr}}}~ \& 
Corentin Le Bihan\footnote{UMPA, ENS Lyon; {\color{blue} {\tt corentin.le-bihan@ens-lyon.fr}} }}

\maketitle
\section{Introduction}

Let $X$ be a compact smooth\footnote{In all this paper, ``smooth'' means $C^\infty$} manifold without boundary equipped with a smooth density $|dx|$.
We will often denote by $L^2$ the Hilbert space $L^2(X,|dx|)$. 
Let $P$ be a differential operator of degree $2$ on $X$ with smooth coefficients acting 
on complex valued functions. We assume in what follows that $P$ is 
symmetric, i.e., for any pair of smooth  complex valued functions $f, g$ on $X$, we have 
$\int _X Pf ~\bar{g} |dx|= \int_X f~ \overline{Pg}|dx|$. 

The adjoint $P^\star $ of $P$ is then defined as follows: the domain $D(P^\star)$ is the set of distributions $f\in L^2$
so that $Pf$ belongs to $L^2$ and $P^\star $ is the operator $P$ acting on such distributions. 
\begin{defi} A symmetric linear differential operator with smooth coefficients
$P$ is {\rm essentially self-adjoint}
 (denoted ESA in what follows) if the graph of $P^\star $ in $L^2\times L^2 $ is the closure of the graph of
$P$ with domain the smooth functions on $X$. More explicitely,  for each $v=Pu $ with $u,v \in L ^2$,
 there exists a sequence $(u_n, v_n)$ with $u_n $ smooth, $v_n =Pu_n$
and $(u_n, v_n )\ra (u,v)$ in $L^2\times L^2$.

\end{defi}
The ESA property is also called {\it Quantum completeness} because the evolution equation
$du/dt =iPu $ with $u(t=0) =f $, with $f$ smooth,  has then an unique solution defined for $t\in\R $ and denoted
$u(t)={\rm exp }(itP) f $ (see \cite{R-S-75}, sec. VIII).

On the other hand, $P$ admits a principal symbol and also a sub-principal symbol:
if one chooses local coordinates $x=(x_1, \cdots, x_n)$ so that 
$|dx|= |dx_1 \cdots dx_n |$, and if 
$P=\sum \frac{\pa }{\pa x_k} a_{kl}(x)\frac{\pa }{\pa x_l} + \sum b_k(x) \frac{\pa }{\pa x_k}+ c (x)$, the principal symbol
is $p_2:=-\sum a_{kl}(x) \xi_k \xi_l $           and the sub-principal symbol is
$p_1:=i \sum b_k(x) \xi_k$ (see \cite{G-K-L-64} sec. 3). Note that, if $P$ is symmetric, $p_1$         and $p_2$ are real valued. 
We denote $p=p_2+p_1$ and call it the symbol of $P$. Note that $p$ is real valued if $P$ is formally symmetric.
The symbol $p$ is independent of the choice of local coordinates as soon as interpreted as a function on the cotangent
space $T^\star X$. 
The cotangent space is a symplectic manifold and one can use $p$ as an Hamiltonian function on it.
We say that $P$ is {\it classically complete} if the Hamiltonian flow of $p$ is complete: it means that the maximal interval
of definition of any  integral curve  of the
Hamiltonian vector field of $p$ is  $\R $.

{\bf A natural question is then: how are classical and quantum completeness related?}
The goal of this article is to give very partial answers to this question.
We will state a possible answer as the 
 
\begin{conj}\label{conj:main}  Let $P$ be a formally self-adjoint  differential operator of degree $2$  on $C^\infty (X)$
 where $X$ is a closed smooth manifold equipped with a smooth density
$|dx|$, classical and quantum completeness are equivalent. 
\end{conj}
As we will see, this  conjecture holds true in the following cases:
\begin{enumerate}
\item Differential operators of degree $2$ on the circle
of the form
 \[ P=a(x)d_x^2 +\cdots  \]
 where all zeroes of $a$ are of finite order.
\item Differential operators of degree $1$.
\item  Generic and conformally flat Lorentzian   Laplacians  on 2-tori. 
\end{enumerate}

Let us describe other known results on this question: 
a classical result is  that classical completeness and quantum completeness are not equivalent in general;
 examples of Schr\"odinger operators
on $\R$ are given in
\cite{R-S-75}, section X-1, pp 155-157. However the potentials involved are quite complicated near infinity: they do not admit a polynomial
asymptotic behaviour. 
A classical results in this domain is Gaffney's Theorem \cite{Ga-54} which states that,
if a Riemannian manifold $(X,g)$ is complete, the Laplace operator on it is ESA. For a clear  proof, see \cite{Dav-89}, pp 151--152. 
For  a more recent work on this aspect, see \cite{BMS-02}. 

{\bf Acknowledgements. }
{\it Many thanks to Christophe Bavard, Yves Carri\`ere, Etienne Ghys, Nicolas Lerner and Bernard Malgrange 
for several discussions  helping us. Many thanks also to the referee for many remarks helping us to make the paper more
accessible. }
\section{General facts  on ESA operators }
\subsection{Abstract context}\label{ss:abstract}
Let us recall some classical  results which  can be found in \cite{R-S-75}, sec. X.
Let  $({\cal H},<.|.>)$ be an Hilbert space. A linear operator $P:D(P)\ra {\cal H}$ with $D(P)$ a dense subspace of ${\cal H}$
is said to be symmetric if, for  all $x,y\in D(P)$, $<Px|y>=<x|Py > $. The closure of $P$ is the operator
$\bar{P}$ whose graph is the closure of the graph of $P$ in ${\cal H}\times {\cal H}$. 
The adjoint $P^\star $ of $P$ is defined as follows: the domain $D(P^\star )$ is the set
of $x\in {\cal H}$ so that $y \ra <Py |x > $ as defined on $D(P)$ extends continously to ${\cal H}$. We have
then $ <Py |x >\equiv <y|z> $ and we define $P^\star x =z$. 
The operator $P$ is ESA if $P^\star $ is the closure of $P$. In other words, $P$ has an unique self-adjoint extension.
A useful property is the following one
\begin{theo} A symmetric operator $P$ on an Hilbert space ${\cal H}$ 
 is ESA if and only if  the spaces $\ker_{\cal H} (P^\star \pm i) $ are $\{0 \}$. 
\end{theo}

\subsection{The case of differential operators on compact manifolds without boundary}
Recall that a differential operator with smooth coefficients acts on Schwartz distributions and in particular on $L^2$
functions. 
In particular, we see that any  symmetric elliptic operator $P$ on a closed manifold
 is ESA: if $(P\pm i )u=0$, $u$ is smooth and the result follows from the
symmetry of $P$. This is why we are only interested here to non elliptic operators. 
\subsection{The case of differential operators on a compact interval }

Let $X:=[\alpha ,\beta]$  be  a compact   interval. 
 We consider  a  differential operator $P$  of degree $2$ whose  coefficients are smooth up to the boundary. 
Assuming that $P$ is symmetric on  $C_c^\infty (]\alpha,\beta[, |dx|)$ (it is usually called formally symmetric), then  $P$ is given by the 
equation (\ref{equ:P}). 

  We want to describe   the {\it Dirichlet boundary
conditions}. For that, we  assume that $P$ is elliptic near the boundary, i.e. that $a$ does not vanish at the points of $\pa X$. 
We will take for domain of $P$ the space $D(P):=C^\infty (\bar {X},\C )\cap \{ f|f|_{\pa X}=0 \}$. 
\begin{lemm}
The domain of $P^\star $ is then the set of $L^2$ functions $g$ so that $Pg \in L^2$, where $P$ is acting on distributions defined
in the interior of $X$  and $g$ vanishes on the boundary. 
\end{lemm}
{\it Proof.--}We get first that $Pg\in L^2$ by  looking at $\int_X Pf \bar{g} |dx| $ with $f\in C_c^\infty (]\alpha, \beta [)$. 
It  follows that $g$ is continuous  near the boundary.
Then we have, if $f\in D(P)$, 
\[ \int  _X (Pf \bar{g} - f \overline{Pg}) |dx| =  [af'g]_{\alpha }^\beta \]
We have to control the righthandside  in terms of the $L^2$ norm of $f$ which is clearly not
possible if $g$ does not vanish on the boundary because $a$ does not vanish at on $\pa X$.  
\hfil $\square$
\subsection{Localization}
Let us prove the following localization 

\begin{lemm}\label{lemm:loc}  Let $P$ be a symmetric  operator of degree $2$ on a the circle. 
 Let $Z\subset X$ be the closed set of points where $P$ is not elliptic.
We assume that $Z$ is a finite set $Z=\{x_1,\cdots,x_N \}$ .   Let $\Omega =\cup_{j=1}^N [\alpha_j ,\beta_j] $ be a neighbourhood of $Z$
so that $P$ is elliptic near the boundary of $\Omega $. 
 Then $P$ is ESA if and only if the Dirichlet restriction $P_\Omega$ of $P$ to $\Omega$ is ESA.
\end{lemm}
{\it Proof.--}

Let us first prove that, if $P_\Omega$ is ESA, $P$ is ESA:
let us take a $\rho \in C_c^\infty (\Omega) $ with $\rho \equiv 1$ near $Z$.
Then, if $Pu=v $ with $u,v \in L^2 (X)$, $(1-\rho) u $ belongs to the Sobolev space
$ H^2(X)$ by ellipticity of $P$ on the support of $1-\rho$. In particular $P((1-\rho)u)\in L^2$. 
There exists $(u'_n,v_n'=Pu_n')$ a sequence of smooth functions converging to $((1-\rho) u ,P((1- \rho) u )$ in $L^2$ by density
of $C^\infty (X)$ in $H^2(X)$. We have now $P(\rho u ) = w $ with $\rho u, w \in L^2 $ and ${\rm support }( \rho u )
\subset \Omega $. ESA of $P_\Omega$ allows to approximate $ ( \rho u, w )$ by smooth functions $(u_n'', v_n'' =P u_n'' )$
 and we can assume
that $u_n''$ vanishes near the boundary because $ \rho u$ does.
Then $(u'_n+u''_n ,v'_n+ v''_n)$ are smooth,  converge to $(u,v)$ in $L^2 \times L^2$ and
$P(u'_n+u''_n )=v'_n+ v''_n$. This allows to conclude that $P$ is ESA. 

Let us now  prove that if $P$ is ESA, $P_\Omega$ is ESA:
let us start with $P_\Omega u=v $ with $(u,v)\in L^2(\Omega)$, $u\in H^2$ near the boundary
 and $u(\pa \Omega)=0$.
We choose $\rho \in C_c^\infty (\Omega )$ with $\rho =1$ near $Z$. 
 Similarly to the previous argument, we decompose
$u= \rho u + (1-\rho)u $. And by ellipticity of $P $ near $\pa \Omega$ we get that $(1-\rho)u $ belongs to the Sobolev
space $H_0^2(\Omega)$ of distributions which are in $H^2(\Omega)$ and vanish at the boundary. By density of smooth functions
vanishing at the boundary  in  the Sobolev space $H_0^2$,
 we get an approximating sequence to $((1-\rho)u , P((1-\rho)u) $.
Now we are left with $\rho u $ with $P(\rho u ) \in L^2$ and we use the fact that $P$ is ESA to get an approximating
sequence $(u'_n, Pu'_n)$ on $X$. Choosing $\rho_1 \in C_c^\infty (\Omega )$  with $\rho_1=1$ on the support of $\rho$,
we take $u''_n =\rho_1 u'_n$ and we get an approximating sequence for $(\rho u, P(\rho u)$. 
This allows to conclude.  
{\hfill  $\square $ }

Note that the previous localization result extends probably to higher dimensional manifolds, but this extension is much less
simple.

\section{Essential self-adjointness of differential operators of degree 1}\label{sec:esa1D}

   The following property goes back to 
 Friedrichs \cite{Fr-44} as cited by H\"ormander
\cite{Hor-65}. 
\begin{lemm} Let $P$ be a differential operator of degree $1$ on a closed manifold $X$  and $u \in L^2(X)$ 
so that $Pu\in L^2(X)$. There exists a sequence $u_j \in C^\infty (X)$ so that
$u_j \ra u $  and $Pu_j \ra Pu $ both in $L^2 (X)$.
\end{lemm}
If $P$ is symmetric, this implies that the closure in $L^2 \oplus L^2$ of the graph of $P $ restricted
to smooth functions is the graph of the adjoint of $P$. Hence $P$ with domain $C^\infty (X)$  is essentially self-adjoint.  

\begin{theo} Any symmetric differential operator of degree $1$ on a closed manifold is
essentially self-adjoint.
\end{theo}
This holds in particular for differential operators of the form
$P:=i(V+ \ha {\rm div}_{|dx|}(V))$ where $V$ is a vector field and 
$ {\rm div}_{|dx|}(V):= d (\iota (V)dx )/ dx $ where  $\iota $ is the inner product. 

 In the note \cite{L-19}, Nicolas Lerner remarks that this property extends to
  pseudo-differential operator of degree $1$. 

This is related to our problem  because then the Hamiltonian flow is complete at infinity: the Hamiltonian vector field of $p$
is bounded by $C\| \xi \|$ and the completeness at infinity follows from Gronwall lemma.  
\section{Sturm-Liouville operators on the circle}\label{sec:sturm}

\subsection{Main result}
Any symmetric operator on the circle  $S^1=\R/\Z$ equipped with the Lebesgue measure $|dx|$ can be
written as 
\begin{equation} \label{equ:P} P =d_x a(x) d_x -ib(x) d_x -i\ha b'(x) + c(x) \end{equation} 
where $d_x :=d/dx $, 
 $a,~b,~c $ are smooth real valued periodic functions of period $1$. We assume always in what follows that 
the zeroes of $a$ are of finite multiplicities.
The symbol $p$ of $P$ is 
\[ p=-a(x)\xi^2+ b(x)\xi  \]
The term $c(x)$ plays no role in the essential self-adjointeness, so we will forget it in what follows.

Our main result is 
\begin{theo}\label{theo:main1D} For operators $P$ of the previous form, classical completeness of the Hamiltonian flow of $p$ is equivalent
to quantum completeness of $P$. \end{theo}
Our proof consists in describing the properties of $a$ and $b$ leading to classical completeness and to study the quantum
completeness in the corresponding cases. 

Note that the result in the case where the zeroes of $a$ are non degenerate is also  proved using some 
microlocal analysis  in \cite{Tai-20}. 

\subsection{Classical completeness}

We have the
\begin{theo}
Let $p:=-a(x)\xi^2 + b(x)\xi $ where the zeroes of $a$ are of finite multiplicity.
  Then the Hamiltonian flow of $p$ is complete on $T^\star S^1$ 
 if and only if the zeroes of $a$ are not simple and 
$b$   vanishes at these zeroes. Moreover, this flow is complete if and only if it is null complete, i.e. complete when
restricted to $p^{-1}(0)$.
\end{theo}
{\it Proof.--} 
Recall that the Hamiltonian differential equation writes 
\[ dx/dt=-2a(x)\xi+ b(x), ~ d\xi/dt =  a'(x)\xi^2- b'(x) \xi .\]
The function $p$ is constant along  the integral curves. 
We will denote by $(x_0,\xi_0)$ the data at time $0$ of the integral curves that we will consider. 

The proof splits into three cases:
 \begin{enumerate}

\item {\it Assume that $a(0)=0 $ and $b(0)> 0$.}
We will show that the flow is not null complete. 
 Let us look at the set $p=0$ near $x=0$. This set is the union of the disjoint
curves
$\{ \xi=0 \}$ and $C:=\{ a(x)\xi - b(x)=0 \}$. The curve $C \cap \{ x >0 \} $ is oriented by the flow so that $x$ is decaying
because then $dx/dt=-b(x)$.  
 Let us start on $C$,  with $x_0 >0 $ small enough and $\xi_0$ so that 
$-a(x_0)\xi_0 +b(x_0)= 0$. Then $-a(x(t))\xi(t)+b(x(t))=0$ for all $t$.  Along $C$, we have   $dx/dt =-b(x)$.
Hence, there exists $t_0 >0$ so that $x(t_0)=0$ and 
$\xi(t_0)=+  \infty $. The flow is not null complete: the maximal integral curve is only defined up to $t_0^-$.

\item {\it Assume that $0$ is a non degenerate zero of $a$ and $b(0)=0$.}
Let us start with $x_0=0$ and $\xi_0\ne 0 $. We have $x(t)= 0$ for all $t$ and
$d\xi/dt = a'(0)\xi^2 - b'(0) \xi $. The solution of this differential equation is not defined for all $t$'s because $a'(0)\ne 0 $.
The flow is not null complete.

\item {\it Assume now that zeroes of  $a$  are  degenerate  and $b$ vanishes on $a^{-1}(0)$.} 
We want to prove that the flow is complete.
We have, by conservation of $p$, for any integral curve, there exists $E$ so that  $-a(x)\xi^2+b(x)\xi \equiv E $.
It follows that $\xi $ stays bounded on any compact interval in $x$  disjoint from $a^{-1}(0)$. We need only to consider what happens
when $x(t)$ comes close to a zero of $a$, says  $x=0$. 

Let us first assume that $x_0=0$, then $x(t)\equiv 0$ for all $t$ and $d\xi /dt = -b'(0)\xi $. The trajectory is complete.

If $x_0\ne 0 $ is close to $0$, we get 
$dx/dt = \pm \sqrt{ -4a(x)E -b(x)^2 }= O(|x|)$.
It follows that $x(t)$ does not reach $0$ in finite time and hence  the integral curve ie defined for all times. 
\end{enumerate}
{\hfill $\square $}

\subsection{Simple zeroes of $a$} \label{sec:NDzero}
We will show in this section that if $a(0)=0$ is a simple zero of $a$ then $P$ is not ESA.

Let  $I:=[-\alpha, \alpha ] $ with no other zeroes of $a$ inside $I$.
 The point $0$ is a regular singular point (see Appendix A)
of the differential equation
$(P-i)u=0 $.

The indicial equation writes
$A r^2 -iB r =0$ with $A:=a'(0),~B=b(0)$.
Hence the solutions of this equation near $0$ writes, for $x>0$, 
$y(x)=f(x) + x_+ ^{iB/A } g(x) $ if $B\ne 0$
and 
$y(x)=f(x) +  g(x)\log x  $ if $B=0$ 
with $f,g$ smooth up to $x=0$ (see Appendix A)  and similarly for $x<0$ with $x_-$ and $\log (-x)$.

Let $y_+ $ be the unique solution of $(P-i)y_+ =0,~ y'_+ (\alpha)=0 , ~y'_+ (\alpha )=1$ on $]0,\alpha ]$, so that $y_+$ satisfies 
the Dirichlet boundary condition at $\alpha $.
And define $y_-$ similarly with $y_- (-\alpha)=0$. If we extend $y_+ $ by zero for $x<0$, we get a Schwartz distribution $Y_+$ and
$(P-i)Y_+$ is supported by the origine. We have 
$PY_+= d_x d_x a Y_+ +d_x [a,d_x] Y_+ -ib d_x Y_+-ib'Y_+/2$. We check  that $d_x aY_+ $ is in $L^2_{\rm loc}$.  So  that 
 $(P-i)Y_+ $ is near $0$  in the Sobolev space $ H^{-1}$, because $d_x aY_+ \in L^2$. 
The derivatives $\delta '(0),  \cdots $ of the Dirac distribution are not in $ H^{-1}$.
We have the same result for $Y_-$.
It follows that
$(P-i)Y_\pm =\mu_\pm \delta (0)$.

Hence there is a non zero linear combination $Y$ of $Y_+ $ and $Y_-$ which satisfies
$(P-i)Y=0$ and the Dirichlet boundary conditions at $\pm \alpha $. This proves that $P_I$ is not ESA and hence $P$ is not ESA
by Lemma \ref{lemm:loc}.

\subsection{Degenerate zeroes where $b(0)$ vanishes}

Let us assume that all zeroes of $a$ are degenerate. Then if $I=]c,d[ $ is an interval between two zeroes of
$a$ and assume $a>0$ on $I$, we will show that
there is an explicit unitary map from $L^2(I,|dx|)$ onto $L^2(\R, |dy |)$ sending $C_c^\infty (I)$ (the set of operator with compact support on $I$) into
$C_c^\infty (\R)$ (set of operator with compact support on $\mathbb{R}$) and sending $P$ to 
$Q= d_y^2 + V$ with an explicit $V$.

{\it  First step: a gauge transform.}
Let us consider $P_S:=e^{-iS}P e^{iS}$ where $S$ is smooth and real valued.
We get, by an easy calculation, 
\[ P_S= d_x a d_x -i (b-2a S' ) d_x - i(b'/2 -a'S'-aS") -aS'^2+bS'   \]
Choosing $S $ so that $S'=b/2a$, we get
\[ P_S = d_x a d_x +b^2/4a \]

{\it  Second step: a change of variable.}

Let us choose $x_0\in I$.
Let us define $y=\phi(x)=\int _{x_0}^x  a^{-\ha} (t)dt $. The map 
$\phi $  is smooth diffeomorphism in  $I$ onto $\R$. 
Let us introduce the unitary transform 
$\Omega :L^2 (I,|dx| ) \ra L^2 (\R, |dy |)$ defined by
$\Omega f(\phi(x)) =a(x)^{1/4} f(x)$.
We compute
$P_\Omega :=\Omega  P_S \Omega ^\star $.
We get
$P_\Omega = d_y^2 + V(y)$ with 
\[ V(y)= \left( \frac{b^2}{2a} + \frac{ a'^2}{16 a} - \frac{a''}{4} \right)(\phi^{-1}(y))   \]
$(a'^2/16a +a''/4)$ is bounded. The derivative of $(b/a^{1/2})\circ\phi^{-1}$ is 
\[\left(b'-\frac{a'b}{2a}\right)\circ\phi^{-1}\]
which is also bounded. So $(b/a^{1/2})\circ\phi^{-1}$ is bounded by $C(1+|y|)$ and we get that $V(y)\leq C (y^2+1)$. 

It follows then from the Farine-Lavis Theorem (Theorem X.38 of \cite{R-S-75}) that $P_\Omega $ is ESA and hence $P$ is ESA
on $C_c^\infty (I)$. It follows that $P$ is  ESA  on $C_c^\infty (S^1\setminus a^{-1}(0))$ and a fortiori on $C^\infty (S^1)$.

\subsection{Degenerate zeroes where $b(0)$ does not vanish}

Finally, we study the case where all zeroes of $a$ are degenerate and $b$ does not vanish at least at one of these, say
$x=0$. 
We will need  the following
\begin{lemm}
Let us choose a smooth function $E $ on $I:=]0,c]$ so that $E'= b/a $. There exists two independent  solutions $u_1$ and $u_2$
of $(P-i)u=0$ on $I$
such that  $u_1$ is smooth up to $0$,
$u_2= u_3 e^{iE } $ with $u_3$ smooth up to $0$.
\end{lemm}
It follows that the functions $a(x) d_x u_j $ are in $L^2$ and 
 that $P$ is not ESA by the same argument than in Section \ref{sec:NDzero}. 

{\it Proof.--} (of Lemma)
We check first the existence of $u_1$ in an elementary way by showing the existence of a 
 full Taylor expansion directly: we  start with  the Ansatz $u_1(x)=1+a_1x + a_2 x^2 +\cdots $.
We get $ b(0) a_1+ (b'(0)/2 -1) =0 $. Hence $a_1$. Then, inductively, we get an expression for $a_k$ as a function of the 
$a_l$ for $l<k$. Applying Malgrange's Theorem 7.1 in \cite{Ma-74}, we get a smooth solution $u_1$ with $u_1(0)=1$. 

 Then we make the Ansatz 
$u_2= u_1v $ and we get the following differential equation for $v$:
\[ \left( d_x + \frac{a'}{a}+ 2\frac{u'_1}{u_1} -i\frac{b}{a} \right )d_x v =0  \]
It follows that, we can choose
\[ d_x v= \frac{1}{a u_1^2}e^{iE } \]
If $k$ is the order of the zero of $a$ at $x=0$, 
we can  choose local coordinates near $0$ so that $E=1/y^{k-1}$, we get
$d_yv= A(y)y^{-k}{\rm exp} (i/y^{k-1})$. 
We can integrate by part and we get
\[ v(y)= {A_0(y)}e^{i/y^{k-1}} -\int_y^1 A_1(y)e^{i/y^{k-1}}\]
with $A_0$ and $A_1$ smooth. We can iterate the integration by part and get
a formal solution $v(x)\equiv (v_0+v_1x +\cdots )e^{iE (x) }$.
Again we can apply  Malgrange's Theorem.
{\hfill $\square $}


\section{Lorentzian Laplacians on surfaces}

\subsection{General facts on Lorentzian tori}
We will consider for $X$ a 2-torus with  a smooth Lorentzian metric $g$.
Recall that a Lorentzian metric on a surface is a smooth non degenerate symmetric 2-form of signature $(1, 1)$. 
 There is, as in the Riemannian case, an associated 
geodesic flow (the Hamiltonian flow of the dual metric), a canonical volume form and a Laplace operator, which is 
an hyperbolic operator. 

{\it The null curves:} a smooth curve $\gamma $  is said to be {\it null } if, at every point $x$ of $\gamma $,
and any tangent vector $V$ to $\gamma $ at the point $x$, we have
$g_x(V,V)=0$, i.e. the tangent spaces to $\gamma $ are isotropic.
Locally the null curves are the leaves of two transverse foliations. This is not always true globally.

A {\it  closed null leaf} $ \gamma $  is a simple closed curve which is  null. 
There is then a neighborhood of $\gamma $ with two null foliations: ${\cal F}_+ $ is close to the tangent space to $\gamma $
and ${\cal F }_- $ is transversal to $\gamma $.
We can then define a {\it Poincar\'e map} $P_\gamma $ as follows. Take a point $x_0$ on $\gamma $
and a germ of leaf of ${\cal F}_- $,  $C$ at $x_0$. 
Then  $P_\gamma $ is a germ of diffeomorphism $(C,x_0 )$ into itself obtained by following the null leaves of ${\cal F}_+$. 
The map $P_\gamma $ is uniquely defined modulo smooth conjugation by germs of diffeomorphisms. 
Note that $P_\gamma $ is orientation preserving because $X$ is orientable. 

\subsection{Examples of non geodesically complete Lorentzian surfaces}
It is  known that Lorentzian metrics on the 2-torus are not always geodesically complete. It is the case for example
for the  Clifton-Pohl torus: 

Let $T$ be the quotient of $\R^2 \setminus 0 $ by the group generated by the homothety of ratio $2$.
On $T$, the Clifton-Pohl Lorentzian metric is
$g:=dx dy /(x^2+y^2)$. The associated Laplacian 
$\Box_g = (x^2 + y^2) \partial^2 /\pa x \pa y $ is formally self-adjoint on
$L^2 (T, |dx dy| /(x^2+y^2))$.

There is also a much simpler example, namely the quotient on $(\R_x^+ \times \R_y, dx dy )$
by the group generated by $(x,y)\ra (2x, y/2) $. The manifold is not closed, but non completeness sits already in a compact region. 

It is known these metric are  not geodesically complete.
What about ESA  of $\Box _g $?


\subsection{Some results}

We will prove a rather general result: 
\begin{theo}\label{theo:lor}
1) If the metric $g$ admits a closed null leaf for which the Poincar\'e section is not tangent to infinite 
order to the identity, then $g$ is not geodesically complete. Under the same assumptions,   $\Box _g $ is not ESA.

2) If $g$ is conformal to a flat metric with a smooth conformal factor on a 2-torus, then  $\Box _g $ is  ESA.
\end{theo}

\begin{rem}In the first case, the proof of  the null-incompleteness of the geodesic flow
 is due to Yves Carri\`ere and Luc Rozoy \cite{C-R-94}.
 We will reprove it.

Note that the conformal class of a Lorentzian metric is determined by the null foliations; hence ESA is a property
of these foliations. 
\end{rem} 
For the proof of Theorem \ref{theo:lor}, we will  need two lemmas:
\begin{lemm}The null-geodesic completeness  is invariant by conformal change. 
\end{lemm}
{\it Proof.--} If $g=e^\phi g_0$, the dual metric satisfy 
$g^\star =e^{-\phi}g^\star_0 $ and hence the geodesic flow  restricted to $g^\star =0$ are conformal with a bounded ratio.
{\hfill $\square $}
\begin{lemm}The ESA property is invariant by conformal change. 
\end{lemm}
{\it Proof.--}If $\Box _g u=v$, we have also
$ \Box _{g_0}u =e^\phi v $ and $e^\phi v $ is in $L^2$ as soon as $v$ is. Hence, if $\Box _{g_0} $ is ESA, there
exists a sequence $(u_n,w_n = \Box _{g_0}u_n )_{n \in \N}$ converging in $L^2$ to $(u, e^\phi v)$ and 
$e^{-\phi} v_n $ converges to $v$. 
{\hfill $\square $}

 This proves part 2 of Theorem \ref{theo:lor}. 

\subsection{Normal forms}

It is well known and due to Sternberg \cite{St-57} that a smooth germ of map $(\R,0)\ra (\R,0)$ whose differential  at the origin is 
in $]0,1[\cup ]1,+\infty [ $  is smoothly 
conjugated to $y \ra \lambda y $ and hence is the time 1 flow of the vector field $\mu y \pa_y$ with $\lambda =e^\mu$.  

A similar result hold for more degenerate diffeomorphisms:
we assume that $g$ admits a closed null-leaf so that the Poincar\'e map  is of the
form $P(y)=y+ y^k R(y) $ where $R (0)\ne 0$ and   $k\geq 2$.
It is proved in \cite{Ta-73} (see Theorem 4)
\begin{theo}\label{theo:tak} Any such map is the flow at time 1 of a vector field
$V= A(y)\pa _y $ with $A\sim A_0 y^k$. 
\end{theo}

Let $\gamma $ be closed null-leaf of $g$ and $U$ a \nghbd ~of $\gamma $ so that we have the two  null foliations
 ${\cal F}_+ $ and ${\cal F}_-$.
We have the:
\begin{theo} Let $\gamma $ be a closed leaf whose Poincar\'e map is $P={\rm Id}+ R $ with $R$ of order $k$. 
There exists  coordinates near $\gamma $ so that  
the metric $g$ is conformal
to
$g_0=dx (dy -a(y) dx )$ with $a(y)\sim a_0y ^k, a_0\ne 0 $. 
\end{theo} 
{\it Proof.--}
 Let us parametrize  the closed leaf $\gamma $ by $x\in \R/\Z$ and extend the coordinate $x$ in some
\nghbd ~$U$  of $\gamma $ so that the null foliation ${\cal F}_-$  is given by $dx=0$. Choose then for $y$ any coordinate in 
$U$ so  that $y=0$ on $\gamma $. 
We introduce  the differential equation $dy/dx = A(x,y)$ associated to the foliation ${\cal F}_+$   close to $\gamma $. Note that
$A(x,0)=0$. 
Let  $\phi_x (y)$ be the flow of this differential equation. The map $y\ra \phi_1 (y)$ is the Poincar\'e map
of $\gamma $.
By Theorem \ref{theo:tak}, we can  choose a vector field $a(y)\pa y$  so that the time one flow is  the same Poincar\'e map;
and denote by $(\phi_0)_x (y) $ this flow.
Let us consider the germ of diffeomorphism near $\gamma $ defined
by \[F: (x,y) \ra \big(x, y'=(\phi_0)_x \circ \phi_x ^{-1}(y) \big).\] The map $F$ sends the integral curves of $dy -bdx $ onto
the integral curves of $dy-a dx $ and is periodic of period $1$ because the time 1 flows are the same.
 Hence the two null foliation are given respectively by $dx=0$
and $dy'-a(y') dx =0$. The Theorem follows.

{\hfill $\square $}

\subsection{Proof of Theorem \ref{theo:lor}, part 1}
 The idea is to use the normal form which, being invariant by translation in $x$, allows a separation 
of variables and hence application of Theorem \ref{theo:main1D}. 

Let us first prove the null incompleteness. Using the normal form and the conformal invariance of null completeness,
we have to study near $y=0$ the Hamiltonian $h=-\eta (a(y)\eta +\xi )$. The function $\xi $ is a constant of the motion.
Let us take initial conditions with $y_0 >0$, $\xi_0 >0 $, and $a(y_0)\eta _0 +\xi_0=0$.
We have, using that $a(y)\eta +\xi_0 $ stays at $0$,  $dy/dt= 2a(y)\eta =-2\xi_0 $. Hence $y(t)$ vanishes for a finite time
$t_0$ and, we have then $\eta (t_0)=\infty $.  Null incompleteness follows. 

The Lorentzian Laplacian associated to $g=dx(dy-a(y) dy)$ is given by
\[ \Box = \pa _y a(y) \pa _y + \pa ^2_{xy} \] 
Let us look at solutions of $\Box u=v $ of the form
$u(x,y)  = e^{2\pi i x}v(y) $ with $v$ compactly supported near $0$. 
We have
\[ \Box u (x,y)= e^{2\pi i x}\left(  \pa _y a(y) \pa _y  + 2i\pi \pa _y \right)v(y) \] 
The operator $P:=\pa _y a(y) \pa _y  + 2i\pi \pa _y$ is a Sturm-Liouville operator already studied in section \ref{sec:esa1D}.
$P$ is not ESA. It follows then that there exists $v $ compactly supported near $0$ and $L^2$
so that $Pv =w \in L^2$ and  there is no sequences $(v_n, w_n =P v_n)$ converging in $L^2 \times L^2$ to
$(v,w)$. The result follows. 

\subsection{Genericity}
The goal of this section is to show that, for a generic Lorentzian metric on the 2-torus, there exists at least 
one null closed curve $\gamma $ whose Poincar\'e map is hyperbolic, i.e. the differential of $P_\gamma $ at the point
$x_0$ of $\gamma $ is not tangent to the identity. It follows that, for a generic metric on the 2-torus,
the geodesic flow is not complete and the Lorentzian Laplacian is not ESA. 

A $C^\infty -$generic property is a property which holds for an open dense subset of the metrics in the $C^\infty $
topology. 

We have the
\begin{prop} The existence of  a closed null hyperbolic curve is a $C^\infty -$generic property  of Lorentzian metrics on 
2-tori. 
\end{prop}
The following argument is due to Etienne Ghys. 

{\it Proof.--}
The openness of the set of metric with a closed null hyperbolic geodesic    is evident.

We
say that {\it the metric $g$ splits} if the null leaves belongs to two distincts foliations ${\cal F}_+$
and ${\cal F}_-$. 
We say that {\it the metric $g$ is orientable} if there is a smooth  non vanishing  vector field $V$ on $X$ so that
$g(V,V) $ is strictly positive everywhere. Any orientable metric splits: the two foliations are the boundaries of the
connected component $C_+$ of the cone $g>0$ containing $V$, indeed  using an orientation of $X$, we choose
${\cal F}_+ $ so that the frame generated by ${\cal F}_+ $ and $V$ is positively oriented.
We now study the two  different cases. 
{\it 1. Case where $g$ splits:}
the genericity then follows from the fact that having an hyperbolic closed leaf is a generic property for
a foliation of a torus (see \cite{PdM-82}), here for ${\cal F}_+$. 

{\it 2. Case where $g$ do not split:}
we introduce in this case  a two-fold cover $Y$  of $X$ for which the lift $G$ of the metric $g$ is orientable. 
 This cover $Y$ is equipped with
an involution $J$ exchanging  the two null foliations.
Let  us take a null closed curve $\gamma $ of one of these foliations. Then $J(\gamma )$ is a null closed curve of
the other. They cannot cross:   they have the same rotation number, because $J $ is homotopic to the identity. Moreover 
 all intersections have the same sign because both foliations as well as $Y$ are orientable.
 It follows that the projection of $\gamma $ onto $X$ is simple.
 Still having a closed hyperbolic leaf for the foliation ${\cal F}_+ $ of $G$ is generic and ${\cal F}_- =J({\cal F}_+)$. 
{\hfill $\square $}
 
\section{Further questions}

There are still several open problems in this setting; we see at least five of them:

\begin{enumerate}
\item Prove our conjecture \ref{conj:main}. 
\item Describe the self-adjoint extensions in the case of Lorentzian tori in a geometrical way.
\item If we choose a self-adjoint extension, are there interesting spectral asymptotics?
\item Extend to higher dimensional Lorentzian manifolds
\item Extend  to pseudo-differential operators of principal type. 

\end{enumerate}

\begin{appendix}

\section{Appendix: Regular singular points of linear differential equations of order two}
For this section, one can look at \cite{Co-Le-55} and \cite{Wa-65}.

We  consider a linear differential equqation 
\[ Pu:=(a (x) d_x^2 + b (x) d_x + c(x))u=0 \]
We assume that $a(0)=0$ and $0$ is a zero of finite order $k$ of $a$.
The singular point $x=0$ of $P$ is said to be {\it regular } if $b $ (resp. $c$) vanishes at order at least $k-1$ (resp.
$k-2$) at $x=0$.  Otherwise $0 $ is an {\it irregular} singular point.

If $x=0$ is a regular singular point, we introduce the {\it indicial equation}:
\[ a^{(k)}(0)r(r-1)+kb^{(k-1)}(0) r +k(k-1)c^{(k-2)},~r\in \C \]
We call $r_1, ~ r_2$ the two roots of the indicial equation.
Then the following holds:
\begin{itemize}
\item If $\mathrm{Im}(r_1- r_2) \notin \Z$, there exist two  independent solutions of $Pu=0$ on a small interval $]0,c[$ of the form
$u_j=x_+^{r_j}v_j (x)$
\item If $\mathrm{Im}(r_2-r_1) \in \N$, we have
$u_1=x_+^{r_1}v_1(x) $
and $u_2= x_+^{r_2}( v_2(x)\log x + v_3(x))$
\end{itemize}
where the functions $v_j$ are smooth on $[0,c[$ and $v_1(0)=v_2(0)=1$.

\end{appendix}

\bibliographystyle{plain}

\end{document}